\documentclass[12pt]{amsart}
\usepackage{amssymb,latexsym}
\usepackage{array}
\usepackage{enumerate}

\newcommand{\lr}{\mbox{$\longrightarrow$}}

\newcommand{\be}{\begin{equation}}
\newcommand{\ee}{\end{equation}}

\newtheorem{guess}{Theorem}[section]
\newcommand{\bth}{\begin{guess}$\!\!\!${\bf .}~}
\newcommand{\eeth}{\end{guess}}
\renewcommand{\bar}{\overline}
\newtheorem{propo}[guess]{Proposition}
\newcommand{\bpropo}{\begin{propo}$\!\!\!${\bf .}~}
\newcommand{\epropo}{\end{propo}}

\newtheorem{lema}[guess]{Lemma}
\newcommand{\blem}{\begin{lema}$\!\!\!${\bf .}~}
\newcommand{\elem}{\end{lema}}

\newtheorem{defe}[guess]{Definition}
\newcommand{\bdefe}{\begin{defe}$\!\!\!${\bf .}~}
\newcommand{\edefe}{\end{defe}}

\newtheorem{coro}[guess]{Corollary}
\newcommand{\bcor}{\begin{coro}$\!\!\!${\bf .}~}
\newcommand{\ecor}{\end{coro}}

\newtheorem{rema}[guess]{Remark}
\newcommand{\brem}{\begin{rema}$\!\!\!${\bf .}~\rm}
\newcommand{\erem}{\end{rema}}

\newtheorem{exam}[guess]{Example}
\newcommand{\beg}{\begin{exam}$\!\!\!${\bf .}~\rm}
\newcommand{\eeg}{\end{exam}}

\newtheorem{notn}[guess]{Notation}
\newcommand{\bnot}{\begin{notn}$\!\!\!${\bf .}~\rm}
\newcommand{\enot}{\end{notn}}

\newcommand{\cf}{{\mathcal F}}

\newcommand{\cR}{{\mathcal R}}

 \newcommand{\ra}{\rightarrow}
 \newcommand{\da}{\Big
  \downarrow} 
\newcommand{\hra}{\mbox{{$\hookrightarrow$}}}

\newcommand{\cl}{{\mathcal L}}

\newcommand{\bz}{\mathbb{Z}} 

\newcommand{\cv}{{\mathcal V}}

\newcommand{\cz}{{\mathcal Z}}

\newcommand{\cc}{{\mathcal C}} 
\newcommand{\ca}{{\mathcal A}}
\newcommand{\cb}{{\mathcal B}} 
\newcommand{\cx}{{\mathcal X}}
\newcommand{\bq}{\mathbb{Q}}

\newcommand{\bn}{\mathbb{N}}
\newcommand{\br}{\mathbb{R}}

\newcommand{\bc}{\mathbb{C}} 

\newcommand{\co}{{\mathcal O}} 
\newcommand{\ba}{\mathbb{A}}

\renewcommand{\phi}{\varphi}
\newcommand{\tT}{{\widetilde{T}}}
\newcommand{\tG}{{\widetilde{G}}}
\newcommand{\tB}{{\widetilde{B}}}
\newcommand{\tC}{{\widetilde{C}}}
\newcommand{\tW}{{\widetilde{W}}}
\newcommand{\tphi}{{\widetilde{\Phi}}}

\begin{document}

\baselineskip=17pt

\title[Equivariant $K$-theory of regular
compactifications]{Equivariant $K$-theory of regular
  compactifications: further developments}

\author{V. Uma} 
\address{Department of Mathematics, IIT Madras,
  Chennai, India} 
\email{vuma@iitm.ac.in} 
\date{}

    

\begin{abstract}
  In this article we describe the $\tG\times \tG$-equivariant $K$-ring
  of $X$, where $\tG$ is a {\it factorial} cover of a connected
  complex reductive algebraic group $G$, and $X$ is a regular
  compactification of $G$. Furthermore, using the description of
  $K_{\tG\times \tG}(X)$, we describe the ordinary $K$-ring $K(X)$ as
  a free module of rank the cardinality of the Weyl group, over the
  $K$-ring of a toric bundle over $G/B$, with fibre the toric variety
  $\bar{T}^{+}$, associated to a smooth subdivision of the positive
  Weyl chamber. This generalizes our previous work on the wonderful
  compactification (see \cite{u}). Further, we give an explicit
  presentation of $K_{\tG\times \tG}(X)$ as well as $K(X)$ as an
  algebra over the $K_{\tG\times \tG}(\bar{G_{ad}})$ and
  $K(\bar{G_{ad}})$ respectively, where $\bar{G_{ad}}$ is the
  wonderful compactification of the adjoint semisimple group
  $G_{ad}$. Finally, we identify the equivariant and ordinary
  Grothendieck ring of $X$ respectively with the corresponding rings
  of a canonical toric bundle over $\bar{G_{ad}}$ with fiber the toric
  variety $\bar{T}^+$.
\end{abstract}

\maketitle

\thispagestyle{empty} 

\frenchspacing

\textwidth=13.5cm
\textheight=23cm
\parindent=16pt
\oddsidemargin=-0.5cm
\evensidemargin=-0.5cm
\topmargin=-0.5cm

\noindent

\section{Introduction}
Let $G$ denote a connected complex reductive algebraic group,
$B\subset G$ a Borel subgroup and $T\subset B$ a maximal torus of
dimension $l$. Let $C$ be the center of $G$ and let $G_{ad}:=G/C$ be
the corresponding semisimple adjoint group. Let $W$ denote the Weyl
group of $(G,T)$.
 
A normal complete variety $X$ is called an {\it equivariant
  compactification} of $G$ if $X$ contains $G$ as an open subvariety
and the action of $G\times G$ on $G$ by left and right multiplication
extends to $X$.  We say that $X$ is a {\it regular compactification}
of $G$ if $X$ is an equivariant compactification of $G$ which is
regular as a $G\times G$-variety in the sense of \cite[Section
3]{bdp}.  Smooth complete toric varieties are exactly the regular
compactifications of the torus. For the adjoint group $G_{ad}$, the
wonderful compactification $\bar {G_{ad}}$ constructed by De Concini
and Procesi in \cite{dp} is the unique regular compactification of
$G_{ad}$ with a unique closed $G_{ad}\times G_{ad}$-orbit.

We now recall some facts and notations from \cite[Section 3.1]{br}.

Let $\bar{T}$ denote the closure of $T$ in $X$. For the left action of
$T$ (i.e. for the action of $T\times \{1\}$), $\bar{T}$ is a smooth
complete toric variety.  Let $\cf$ be the fan associated to $\bar T$
in $X_{*}(T)\otimes \br$. Then $\cf$ is a smooth subdivision of the
fan associated to the Weyl chambers in $X_{*}(T)\otimes \br$.
Moreover, $\cf=W\cf_{+}$ where $\cf_{+}$ is the subdivision of the
positive Weyl chamber formed by the cones in $\cf$ contained in this
chamber.  Let $\bar{T}^{+}$ denote the smooth toric variety associated to the
fan ${\cf}_{+}$ (for details see \cite[Proposition A1 and Proposition A2]{br}).

Recall that there exists an exact sequence \be\label{fc} 1\ra \cz\ra
\tG:=\tC \times G^{ss}{\stackrel{\pi}{\lr}}G \ra 1\ee where $\cz$ is a
finite central subgroup, $\tC$ is a torus and $G^{ss}$ is semisimple
and simply-connected.  In particular, $\tG$ is {\it factorial} and
$\tB:=\pi^{-1}(B)$ and $\tT:=\pi^{-1}(T)$ are respectively a Borel
subgroup and a maximal torus of $\tG$. (see \cite[Corollary 3.7]{iv}
and \cite{mer}).

In this article we consider the ${\tG\times \tG}$-equivariant $K$-ring
of $X$. We take the natural action of $\tG\times \tG$ on $X$ through
the canonical surjection to $G\times G$. Here we mention that the main
purpose of going to the factorial cover $\tG$ of $G$, is in order to
apply the description of the $\tG\times \tG$-equivariant $K$-ring to
describe the ordinary $K$-ring of $X$. Indeed, the relation of
$K_{\tG\times \tG}(X)$ to ordinary $K$-theory is simpler than that of
$K_{G\times G}(X)$, in view of the results of \cite{mer}.

Recall that, in \cite[Section 3]{u}, taking $G$ to be the simply
connected cover of $G_{ad}$ and $T$ a maximal torus of $G$, we
describe the $G\times G$-equivariant $K$-ring of
$\overline{G_{ad}}$. More precisely, we give $K_{G\times
  G}(\overline{G_{ad}})$ an $R(T)\otimes R(G)$-algebra structure and
determine the multiplication rule explicitly (see \cite[Theorem
3.8]{u}). We further apply this structure to describe the ordinary
Grothendieck ring $K(X)$ as a $K(G/B)$-algebra (see \cite[Theorem
3.12]{u}).

Our primary aim in this article is to generalize the above mentioned
results on the wonderful compactification to any regular
compactification $X$.  With this aim in view, we first give a
description of $K_{\tG\times \tG}(X)$ as a free module of rank $|W|$
over a canonical subring isomorphic to $K_{\tT}(\bar{T}^{+})\otimes
R(\tG)$ (see Theorem \ref{kdec}). Here $K_{\tT}(\bar{T}^{+})$ is the
$\tT$-equivariant $K$-ring of the toric variety $\bar{T}^{+}$, which
is further identified with the subring of $K_{\tG\times \tG}(X)$
generated by $Pic_{\tG\times \tG}(X)$, as an $R(\tG)\otimes
R(\tG)$-algebra. Furthermore, we determine a canonical basis for
$K_{\tG\times \tG}(X)$ as $K_{\tT}(\bar{T}^{+})\otimes R(\tG)$-module,
and describe the multiplication rule of the basis elements (see
Theorem \ref{modkdec}). In particular, when $X$ is the wonderful
compactification of $G_{ad}$, the toric variety $\bar{T_{ad}}^{+}$ is the
affine space $\ba^{r}$, so that $K_{\tT}(\bar{T_{ad}}^{+})$ can be
identified with $R(\tT)$. Thus we get back \cite[Theorem 3.8]{u} in
this case. Note that in the current notation, $\tG$ denotes the simply
connected cover of $G_{ad}$ and $\tT \subseteq\tG$ is a maximal torus.

We bring to notice here that \cite[Section 2]{u} contains general
results on $G\times G$-equivariant $K$-ring of a regular embedding,
which are closely related to the results in this paper.  In
particular, Proposition \ref{key} in this article is an extension of
\cite [Proposition 2.5]{u} to the setting of $\tG\times
\tG$-equivariant $K$-theory (see \cite[Remark 2.7]{u}). Moreover, in
Proposition \ref{key}, we in fact show that we have the inclusion
$K_{\tG\times \tG}(X) \subseteq K_{\tT}(\bar{T}^{+})\otimes R(\tT)$,
as an $R(\tT)\otimes R(\tG)$-algebra. This is sharper than the
inclusion $K_{\tG\times \tG}(X)\subseteq R(\tT)^{\cf_{+}(l)}\otimes
R(\tT)$, analogous to that in \cite [Proposition 2.5]{u}. Indeed, the
inclusion $K_{\tG\times \tG}(X) \subseteq K_{\tT}(\bar{T}^{+})\otimes
R(\tT)$, is key in proving Theorem \ref{kdec}, that generalizes
\cite[Lemma 3.2]{u} and \cite[Theorem 3.3]{u} on the wonderful
compactification, to any regular embedding.

Using the above structure of $K_{\tG\times \tG}(X)$ and applying
(\cite[Theorem 4.2]{mer}), in Section 3, we describe the ordinary
$K$-ring of $X$. More precisely, we prove in this section that the
subring generated by $Pic(X)$ in $K(X)$ is canonically isomorphic to
$\cR(\bar{T}^{+}):=\bz\otimes_{R(\tG)}K_{\tT}(\bar{T}^{+})$. Further,
in Theorem \ref{main}, we show that $K(X)$ is a free module of rank
$|W|$ over $\cR(\bar{T}^{+})$. We finally construct an explicit basis
of $K(X)$ over $\cR(\bar{T}^{+})$ and also determine the
multiplicative structure constants with respect to this basis.  It is
worth noticing here that
$\cR(\bar{T}^{+})=\bz\otimes_{R(\tG)}K_{\tG}(\tG\times^{\tT}\bar{T}^{+})=K(\tG\times^{\tT}\bar{T}^{+})=K(G\times^B
\bar{T}^{+})$ where $B$ acts on $\bar{T}^{+}$ via its quotient $T$. In
other words, $\cR(\bar{T}^{+})$ is the Grothendieck ring of the toric
bundle over $G/B$ associated with $\bar{T}^{+}$. In particular, when
$X$ is the wonderful compactification of $G_{ad}$,
$\cR(\bar{T_{ad}}^{+})$ can be identified with $K(G/B)$. Thus we get
back \cite[Theorem 3.12]{u} in this case.

In Section 3.1 we compare the equivariant and ordinary $K$-ring of any
regular embedding with that of the wonderful compactification, using
the canonical surjective morphism $f:X\ra\bar{G_{ad}}$. We prove that
$K_{\tG\times\tG}(X)$ is a tensor product of $K_{\tG\times
  \tG}(\bar{G_{ad}})$ (where the $\tG$ action on $\bar{G_{ad}}$ is via
its quotient $G^{ss}$) and the subalgebra generated by $Pic_{\tG\times
  \tG}(X)$, over the subalgebra of $K_{\tG\times \tG}(\bar{G_{ad}})$
generated by $Pic_{\tG\times \tG}(\bar{G_{ad}})$. The corresponding
statement for the ordinary $K$-ring is analogous. Indeed, we give an
explicit presentation for $K_{\tG\times\tG}(X)$ as a $K_{\tG\times
  \tG}(\bar{G_{ad}})$-algebra and $K(X)$ as a
$K(\bar{G_{ad}})$-algebra see Theorem \ref{relwond} and Theorem
\ref{relwond1}.

Finally, Section 4 is an appendix on the equivariant $K$-theory of
toric bundles with fibre a smooth semi-projective toric variety. We
take the left action of a connected complex reductive algebraic group
$G$ on the total space and the base, and assume the bundle projection
to be $G$-invariant. In Theorem \ref{ektb} we give a presentation of
the $G$-equivariant Grothendieck ring of the toric bundle as an
algebra over the $G$-equivariant Grothendieck ring of the base.

This result is applied to give a nice interpretation of the
presentations of $K_{\tG\times\tG}(X)$ as a $K_{\tG\times
  \tG}(\bar{G_{ad}})$-algebra and $K(X)$ as a
$K(\bar{G_{ad}})$-algebra obtained in Theorem \ref{relwond} and
Theorem \ref{relwond1} of \S3.1. Indeed it follows by Theorem
\ref{ektb} that $K_{\tG\times\tG}(X)$ as a $K_{\tG\times
  \tG}(\bar{G_{ad}})$-algebra and $K(X)$ as a
$K(\bar{G_{ad}})$-algebra are respectively isomorphic to the
$\tG\times \tG$-equivariant and ordinary Grothendieck ring of a
canonical toric bundle over $\bar{G_{ad}}$ with fiber the toric
variety $\bar{T}^+$ (see Corollary \ref{eqtorbund}).

{\bf Acknowledgement:} I am grateful to Prof.Michel Brion for a
careful reading and several invaluable comments and suggestions on the
earlier versions of this manuscript.

\subsection{Preliminaries on K-theory}
Let $X$ be a smooth projective complex $G$-variety.  Let $K_{G}(X)$
and $K_{T}(X)$ denote the Grothendieck groups of $G$ and
$T$-equivariant coherent sheaves on $X$ respectively. Recall that
$K_{T}(pt)=R(T)$ and $K_{G}(pt)=R(G)$ where $R(T)$ and $R(G)$ denote
respectively the Grothendieck group of complex representations of $T$
and $G$. The Grothendieck group of equivariant coherent sheaves can be
identified with the Grothendieck ring of equivariant vector bundles on
$X$. Further, the structure morphism $X\ra {Spec}~{\bc}$ induces
canonical $R(G)$ and $R(T)$-module structures on $K_{G}(X)$ and
$K_{T}(X)$ respectively (see \cite[Example 2.1]{mer}).

Let $W$ denote the Weyl group and $\Phi$ denote the root system of
$(G,T)$. We have the subset $\Phi^{+}$ of positive roots and its
subset $\Delta=\{\alpha_1,\ldots, \alpha_r\}$ of simple roots where $r$
is the semisimple rank of $G$. For $\alpha\in\Delta$ we denote by
$s_{\alpha}$ the corresponding simple reflection. For any subset
$I\subset \Delta$, let $W_{I}$ denote the subgroup of $W$ generated by
all $s_{\alpha}$ for $\alpha\in I$. At the extremes we have
$W_{\emptyset}=\{1\}$ and $W_{\Delta}=W$.

Let $\Lambda:=X^*(T)$. Then $R(T)$ (the representation ring of the
torus $T$) is isomorphic to the group algebra $\bz[{\Lambda}]$. Let
$e^{\lambda}$ denote the element of $\bz[{\Lambda}]=R(T)$
corresponding to a weight $\lambda\in {\Lambda}$.  Then
$(e^{\lambda})_{\lambda\in {\Lambda}}$ is a basis of the $\bz$ module
$\bz[{\Lambda}]$. Further, since $W$ acts on $X^*(T)$, on
$\bz[{\Lambda}]$ we have the following natural action of $W$ given by
: $w(e^{\lambda})=e^{w(\lambda)}$ for each $w\in W$ and $\lambda\in
{\Lambda}$. Recall that we can identify $R(G)$ with $R(T)^{W}$ via
restriction to $T$, where $R(T)^W$ denotes the subring of $R(T)$
invariant under the action of $W$ (see \cite[Example 1.19]{mer}).

Now, from (\ref{fc}) it follows that $\tB:=\pi^{-1}(B)$ and
$\tT:=\pi^{-1}(T)$ are respectively a Borel subgroup and a maximal
torus of $\tG$. Further, by restricting the map $\pi$ to $\tT$ we get
the following exact sequence: \be\label{fc1} 1\ra \cz\ra \tT\ra T \ra
1. \ee Let $\tW$ and $\tphi$ denote respectively the Weyl group and
the root system of $(\tG,\tT)$.  Then by (\ref{fc}), it also follows
that $\tW=W$ and $\tphi=\Phi$.  Further we have \be\label{e1} R(\tG)=
R(\tC)\otimes R(G^{ss})\ee and \be\label{e2} R(\tT)\simeq R(\tC)\otimes
R(T^{ss})\ee where $T^{ss}$ is the maximal torus $\tT\cap G^{ss}$ of
$G^{ss}$.

Recall we can identify $R(\tG)$ with $R(\tT)^{W}$ via restriction to
$\tT$, and further $R(\tT)$ is a free $R(\tG)$ module of rank $|W|$
(see \cite[Theorem 2.2]{st}). Moreover, since $G^{ss}$ is semi-simple and simply
connected, $R(G^{ss})\simeq \bz[x_1,\ldots, x_r]$ is a polynomial ring
on the ¨fundamental representations¨ (\cite[Example 1.20]{mer}). Hence
$R(\tG)=R(\tC)\otimes R(G^{ss})$ is the tensor product of a polynomial
ring and a Laurent polynomial ring, and hence a regular ring of
dimension $r+{dim}(\tC)={rank}(G)$ where $r$ is the rank of $G^{ss}$.

We shall consider the $\tT$ and $\tG$-equivariant $K$-theory of $X$
where we take the natural actions of $\tT$ and $\tG$ on $X$ through
the canonical surjections to $T$ and $G$ respectively.

We consider $\bz$ as an $R(\tG)$-module by the augmentation map
$\epsilon:R(\tG)\ra \bz$ which maps any $\tG$-representation $V$ to
${dim}(V)$. Moreover, we have the natural restriction
homomorphisms $K_{\tG}(X)\ra$ $K_{\tT}(X)$ and $K_{\tG}(X)\ra K(X)$
where $K(X)$ denotes the ordinary Grothendieck ring of algebraic
vector bundles on $X$.  We then have the following isomorphisms (see
\cite[Proposition 4.1 and Theorem 4.2]{mer}).

\be\label{eqk1} R(\tT)\otimes_{R(\tG)}K_{\tG}(X)\simeq K_{\tT}(X),\ee

 \be\label{eqk2}  K_{\tG}(X)\simeq K_{\tT}(X)^{W},\ee

\be\label{eqk3} \bz\otimes_{R(\tG)} K_{\tG}(X)\simeq K(X).\ee

Let $R(\tT)^{W_I}$ denote the invariant subring of the ring $R(\tT)$
under the action of the subgroup $W_I$ of $W$ for every $I\subset
\Delta$. Thus in particular we have, $R(\tT)^{W}=R(\tG)$ and
$R(\tT)^{\{1\}}=R(\tT)$. Further, for every $I\subset \Delta$,
$R(\tT)^{W_I}$ is a free module over $R(\tG)=R(\tT)^{W}$ of rank
$|W/W_{I}|$ (see \cite[Theorem 2.2] {st}). Indeed, \cite[Theorem
2.2]{st} which we apply here holds for $R(T^{ss})$.  However, since
$W$ acts trivially on the central torus $\tC$ and hence trivially on
$R(\tC)$ we have \be\label{e3}R(\tT)^{W_I}=R(\tC)\otimes
R(T^{ss})^{W_{I}}\ee for every $I\subseteq \Delta$, and hence we
obtain the analogous statement for $R(\tT)$.

Let $W^{I}$ denote the set of minimal length coset representatives of
the parabolic subgroup $W_I$ for every $I\subset \Delta$. Then
$$W^{I}:=\{w\in W \mid l(wv)=l(w)+l(v) ~\forall~ v\in W_I\}=\{w\in
W\mid w(\Phi_{I}^{+})\subset \Phi^{+}\}$$ where $\Phi_I$ is the root
system associated to $W_I$, where $I$ is the set of simple roots.
Recall (see \cite[p. 19]{h}) that we also have: 
$$W^{I}=\{w\in W\mid l(ws)>l(w)~for~all~s\in I\}.$$

Note that $J\subseteq I$ implies that $W^{\Delta\setminus
J}\subseteq W^{\Delta\setminus I}$.  Let
\be\label{1} C^{I}:=W^{\Delta\setminus I}\setminus (\bigcup_{J\subsetneq
I}W^{\Delta \setminus J}).\ee

Let $\alpha_1,\ldots,\alpha_r$ be an ordering of the set $\Delta$ of
simple roots and $\omega_1,\ldots,\omega_r$ denote respectively the
corresponding fundamental weights for the root system of
$(G^{ss},T^{ss})$.  Since $G^{ss}$ is simply connected, the
fundamental weights form a basis for $X^*(T^{ss})$ and hence for every
$\lambda\in X^*(T^{ss})$, $e^{\lambda}\in R(T^{ss})$ is a Laurent
monomial in the elements $e^{\omega_i}:1\leq i\leq r$.

In \cite[Theorem 2.2]{st} Steinberg has defined a basis
$\{f_{v}^{^I}: v\in W^{I}\}$ of $R(T^{ss})^{W_{I}}$ as an
$R(T^{ss})^W$-module.  We recall here this definition: For $v\in
W^{I}$ let \be\label{st1} p_{v}:=\prod_{v^{-1}\alpha_{i}<0} e^{\omega_{i}}\in
R(\tT).\ee Then \be\label{st2} f_v^{^{I}}:=\sum_{x\in W_{I}(v)\big{\backslash}
  W_{I}} x^{-1}v^{-1}p_{v}\ee where $W_{I}(v)$ denotes the stabilizer
of $v^{-1}p_v$ in $W_{I}$.

We shall also denote by $\{f_v^{^I}:v \in W^{I}\}$ the corresponding
basis of $R(\tT)^{W_{I}}$ as an $R(\tT)^{W}$-module where it is
understood that \be\label{notation} f_v^{^I}:=1\otimes f_v^{^I} \in
R(\tC)\otimes R(T^{ss})^{W_{I}}.\ee

\bnot\label{aproposds} Whenever $v\in C^{I}$ we denote
$f_v^{^{\Delta\setminus I}}$ simply by $f_{v}$. We can drop the
superscript in the notation without any ambiguity since $\{C^{I}:
I\subseteq \Delta\}$ are disjoint. Therefore with the modified
notation \cite[Lemma 1.10]{u} implies that: $\{f_{v}: ~v\in
W^{\Delta\setminus I}=\bigsqcup_{J\subseteq I} C^J\}$ form an
$R(\tT)^{W}$-basis for $R(\tT)^{W_{\Delta\setminus I}}$ for every
$I\subseteq \Delta$. Further, let \be\label{2}
R(\tT)_{I}:=\bigoplus_{v\in C^I}R(\tT)^{W}\cdot f_v .\ee \enot

In $R(T)$ let \be\label{mst}f_{v}\cdot f_{v^{\prime}}=\sum_{J\subseteq
  (I\cup I^{\prime})} \sum_{w\in C^{J}}a^{w}_{v,v^{\prime}}\cdot f_{w}
\ee for certain elements $a^{w}_{v,v^{\prime}}\in R(G)=R(T)^{W}$
$~\forall~ v\in C^{I}$, $v^{\prime}\in C^{I^{\prime}}$ and $w\in
C^{J}$, $J\subseteq (I\cup I^{\prime})$.

\section{Equivariant $K$-theory of regular compactifications}

In this section let $X$ be a projective regular compactification of
$G$. We follow the notations in Section 1 together
with the following:

Let ${\cf}_{+}(l)$ denote the set of maximal cones in the fan
$\cf_{+}$. Then we know that ${\cf}_{+}(l)$ parametrizes the closed
$G\times G$-orbits in $X$. By Prop. A1 of \cite{br}, $X^{T\times T}$
is contained in the union $X_c$ of all closed $G\times G$-orbits in
$X$; moreover all such orbits are isomorphic to $G/B^{-}\times G/B$.
Hence $X^{T\times T}$ is parametrized by ${\cf}_{+}(l)\times W\times
W$.

Note that $\cf_{+}$ is combinatorially complete (see \cite[Definition
6.3]{vv}) since $T$ acts on $\bar{T}^{+}$ with enough limits (see
\cite[Remark 2.4]{u}).  Thus by from \cite[Theorem 6.7]{vv} it follows
that $K_{\tT}(\bar{T}^{+})$ is a projective $R(\tT)$-module of rank
$|\cf_{+}(l)|$. This further implies by \cite[Theorem 1.1]{gu}, that
$K_{\tT}(\bar{T}^{+})$ is a free module of rank $|\cf_{+}(l)|$ over
$R(\tT)$.

\bpropo\label{key} There is a chain of injective morphisms of
$R(\tG)\otimes R(\tG)$-algebras as given below: \be\label{chinc}
K_{\tT}(\bar{T}^{+})\otimes R(\tG)\subseteq K_{\tG\times \tG}(X)
\subseteq K_{\tT}(\bar{T}^{+})\otimes R(\tT)\subseteq
\prod_{\sigma\in\cf_{+}(l)}R(\tT_{\sigma})\otimes R(\tT).\ee Moreover,
the $R(\tT)$-algebra structure of $K_{\tT}(\bar{T}^{+})$, induces an
$R(\tT)\otimes R(\tG)$-algebra structure on $K_{\tG\times
  \tG}(X)$. \epropo {\bf Proof:} Note that we have a split exact
sequence \be\label{ses} 1\ra diag(\tT)\ra \tT\times \tT\ra \tT\ra
1,\ee where the second map is $(t_1,t_2)\ra t_1 t_2^{-1}$, and the
splitting is given by $t\ra (t,1)$.

Recall from \cite[Theorem 2.1]{u} that via restrictions to the
$\tT\times \tT$-fixed points $(w,w')\cdot z_{\sigma}$ for $w,w'\in W$,
$K_{\tT\times \tT}(X)$ consists of all families
$(f_{\sigma})(\sigma\in \cf_{+})$ in $R(\tT)\times R(\tT)$, such that

{\it
\begin{enumerate}

\item[(i)] $f_{\sigma,ws_{\alpha},w's_{\alpha}}\equiv f_{\sigma,
w,w'}\pmod {(1-e^{-w(\alpha)}\otimes e^{-w'(\alpha)})}$ whenever
$\alpha\in\Delta$ and the cone $\sigma\in{\cf}_{+}(l)$ has a facet
orthogonal to $\alpha$, and that

\item[(ii)] $f_{\sigma,w,w'}\equiv f_{\sigma^{\prime}, w,w'}\pmod
{(1-e^{-\chi})}$ whenever $\chi\in X^{*}(T)$ and the cones $\sigma$
and $\sigma^{\prime}\in{\cf}_{+}(l)$ have a common facet orthogonal to
$\chi$. 
 
\end{enumerate}
(In $(ii)$, $\chi$ is viewed as a character of $T\times T$ which is
trivial on $diag(T)$ and hence is a character of $T$.)}

Furthermore, by taking $W\times W$-invariants, \cite[Corollary 2.2]{u}
and (\ref{ses}) above imply that, $K_{\tG\times \tG}(X)$ consists in
all families $(f_{\sigma})(\sigma\in\cf_{+}(l))$ of elements of $
R({\tT\times\{1\}})\otimes R(diag(\tT))$ such that
\begin{enumerate}
\item[$(i)$] $(1,s_\alpha)f_{\sigma}(u,v)\equiv f_{\sigma}(u,v)\pmod
{(1-e^{-\alpha(u)})}$ whenever $\alpha\in\Delta$ and the cone
$\sigma\in\cf_{+}(l)$ has a facet orthogonal to $\alpha$.

\item[$(ii)$] $f_{\sigma}\equiv f_{\sigma^{\prime}}\pmod
{(1-e^{-\chi(u)})}$ whenever $\chi\in X^{*}(\tT)$ and the cones $\sigma$
and $\sigma^{\prime}\in{\cf}_{+}(l)$ have a common facet orthogonal to
$\chi$,
\end{enumerate} 
where $u$ and $v$ denote respectively the variables corresponding to
$R(\tT\times\{1\})$ and $R(diag(\tT))$ respectively. Further, recall
from \cite[Remark 2.4]{u}, that we have a canonical restriction
homomorphism \be\label{rm} K_{\tT}(\bar{T}^{+})\ra
K_{\tT}((\bar{T}^{+})^{T})=\prod_{\sigma\in\cf_{+}(l)}R(\tT_{\sigma})\ee where
$\tT_{\sigma}$ denotes the stabilizer at $z_{\sigma}$ and
$R(\tT_{\sigma})=K_{\tT}(z_{\sigma})$.  Moreover, (\ref{rm}) is
injective and the image of $K_{\tT}(\bar{T}^{+})$ can be identified with the
subset of $\prod_{\sigma\in\cf_{+}(l)}R(\tT_{\sigma})$ consisting of
elements $(a_{\sigma})$ such that $a_{\sigma}\equiv
a_{\sigma^{\prime}}\pmod {(1-e^{-\chi})}$ whenever $\chi\in X^{*}(T)$
and the cones $\sigma , \sigma^{\prime}\in \cf_{+}(l)$ have a common
facet orthogonal to $\chi$. Furthermore, (\ref{rm}) is compatible with
the canonical $R(\tT)$-algebra structure on $K_{\tT}(\bar{T}^{+})$, and the
$R(\tT)$-algebra structure on
$\prod_{\sigma\in\cf_{+}(l)}R(\tT_{\sigma})$ given by the diagonal map
(also see proof of \cite[Proposition 2.5]{u}).

Thus $K_{\tT}(\bar{T}^{+})\otimes R(\tT)$ can be identified with the
$R(\tT)\otimes R(\tG)$-subalgebra of $\prod_{\sigma\in\cf_{+}(l)}
R(\tT_{\sigma})\otimes R(\tT)$, generated by the elements
$(a_{\sigma})\otimes b$ such that
$a_{\sigma}-a_{\sigma^{\prime}}\equiv 0 \pmod{(1- e^{-\chi})}$,
whenever $\sigma$ and $\sigma^{\prime}$ share a facet orthogonal to
$\chi\in X^*(T)$.

Further, since $R(\tG)=R(\tT)^{W}$, from the relations $(i)$ and $(ii)$
it follows that $K_\tT(\bar{T}^{+})\otimes R(\tG)$ is a subring of $K_{\tG\times
  \tG}(X)$. In particular, $K_{\tG\times \tG}(X)$ is an algebra over
$K_\tT(\bar{T}^{+})\otimes R(\tG)$, and hence over $R(\tT)\otimes R(\tG)$.

Furthermore, from the relation $(ii)$ above, any element
$$f(u,v)\in \prod_{\sigma\in\cf_{+}(l)}R(\tT_{\sigma})\otimes R(\tT)$$
belonging to $K_{\tG\times \tG}(X)$, must in fact lie in
$K_{\tT}(\bar{T}^{+})\otimes R(\tT)$.  Hence the inclusion $K_{\tG\times
  \tG}(X)\subseteq K_{\tT}(\bar{T}^{+})\otimes R(\tT)$ of $R(\tT)\otimes
R(\tG)$-algebras.$\Box$

\bth\label{kdec} The ring $K_{\tG\times \tG}(X)$ has the following
direct sum decomposition as $K_\tT(\bar{T}^{+})\otimes R(\tG)$-module: \be\label{dsd}
K_{\tG\times \tG}(X)=\bigoplus_{I\subseteq \Delta}\prod_{\alpha\in I}
(1-e^{\alpha(u)})\cdot K_{\tT}(\bar{T}^{+})\otimes R(\tT)_{I}.\ee Further, the
above direct sum is a free $K_\tT(\bar{T}^{+})\otimes R(\tG)$-module of rank
$|W|$ with basis $$\{\prod_{\alpha\in I}(1-e^{\alpha(u)})\otimes f_{v}
:~v\in C^{I}~and~ I\subseteq \Delta\},$$ where $C^{I}$ is as defined
in (\ref{1}) and $\{f_{v}\}$ is as defined above. Moreover, we can
identify the component $K_{\tT}(\bar{T}^{+})\otimes 1\subseteq K_\tT(\bar{T}^{+})\otimes
R(\tT)^{W}$ in the above direct sum with the subring of $K_{\tG\times
  \tG}(X)$ generated by generated by $Pic^{\tG\times \tG}(X)$. 
\eeth

{\it Proof :} Recall from Notation \ref{aproposds} that we have the
following
 decompositions as
$R(\tT)^W$-modules: \be\label{dec1} R(\tT)=\bigoplus_{I}R(\tT)_{I}\ee
\be\label{dec2} R(\tT)^{W_{\Delta\setminus I}}=\bigoplus_{J\subset
  I}R(\tT)_{J}\ee Let
\be\label{fds} L:=\bigoplus_{I\subset \Delta}\prod_{\alpha\in I}(1-e^{\alpha(u)})
K_{\tT}(\bar{T}^{+})\otimes R(\tT)_{I}.\ee For $I\subseteq \Delta$, the piece
$\prod_{\alpha\in I}(1-e^{\alpha(u)}) K_{\tT}(\bar{T}^{+})\otimes R(\tT)_{I}$ in
the above direct sum decomposition is a free $K_\tT(\bar{T}^{+})\otimes
R(\tG)$-module with basis $\{1\otimes f^{I}_{v}~:~v\in C^{I}\}$, 
$f^{I}_{v}$ as in (\ref{st2}). In particular, $L$ is a free
$K_{\tT}(\bar{T}^{+})\otimes R(\tG)$-module of rank $|W|$, and hence is free
over $R(\tT)\otimes R(\tG)$ of rank $|W|\cdot |\cf_{+}(l)|$.

We observe that an element in $\prod_{\alpha\in I}(1-e^{\alpha(u)})
K_{\tT}(\bar{T}^{+})\otimes R(\tT)_{I}$, clearly satisfies the relations $(i)$
and $(ii)$, which define $K_{\tG\times \tG}(X)$ as a subalgebra of
$\prod_{\sigma\in\cf_{+}(l)} R(\tT_{\sigma})\otimes R(\tT)$. We
therefore have an inclusion $L\subseteq K$ of $K_{\tT}(\bar{T}^{+})\otimes
R(\tG)$-modules, where $K:=K_{\tG\times \tG}(X)$.  Further, since
$K_{\tT}(\bar{T}^{+})$ is an $R(\tT)$-algebra $L\subseteq K$ is also an
inclusion as $R(\tT)\otimes R(\tG)$-modules.

Now, by \cite[Lemma 1.6]{u}, it follows that $K$ is a free module over
$R(\tG)\otimes R(\tG)$ of rank $|W|^2\cdot |\cf_{+}(l)|$. Moreover,
since $K$ is a $R(\tT)\otimes R(\tG)$-algebra by Proposition
\ref{key}, it follows that $K$ is a projective module over
$R(\tT)\otimes R(\tG)$. (In view of \cite[Theorem 1.1]{gu}, and using
the fact that $R(\tT)\otimes R(G)$ is free over $R(\tG)\otimes R(\tG)$
of rank $|W|$, this further implies that $K$ is free over
$R(\tT)\otimes R(\tG)$ of rank $|W|\cdot |\cf_{+}(l)|$.)

Thus, $L\hookrightarrow K\ra K/L\ra 0$ is a short exact sequence of
$K_{\tT}(\bar{T}^{+})\otimes R(\tG)$ modules. Moreover, since $K$ and
$L$ are projective as $R(\tT)\otimes R(\tG)$-modules, it follows that
$K/L$ is of projective dimension $1$ as a module over $R(\tT)\otimes
R(\tG)$.

We now require to prove that $L=K$. This we do in the
following lemma, using the inclusion $L\subseteq K$ as $R(\tT)\otimes
R(\tG)$-modules.

\blem{\label{loc}} Let $t_{\alpha}:=\prod_{\beta\neq \alpha}
(1-e^{-{\beta(u)}})\in R(\tT)\otimes R(\tG)$ for every
$\alpha\in\Delta$. Then, the localization $(K/L)_{t_{\alpha}}=0$ for every
$\alpha\in\Delta$.  \elem {\it Proof :} Let \be\label{ab1} M_{\alpha}:=K_T(\bar{T}^{+})\otimes
R(T)^{s_{\alpha}} \bigoplus (1-e^{-{\alpha(u)}}) K_{\tT}(\bar{T}^{+})\otimes
e^{\omega_{\alpha}(v)}.R(\tT)^{s_{\alpha}}.\ee
Further, $R(\tT)=R(\tT)^{s_\alpha}\bigoplus
e^{\omega_{\alpha}}R(\tT)^{s_{\alpha}}$ for every $\alpha\in\Delta$
where $\omega_{\alpha}$ is the fundamental weight corresponding to
$\alpha\in\Delta$. Hence \be\label{ab2} K_{\tT}(\bar{T}^{+})\otimes R(\tT)= K_\tT(\bar{T}^{+})\otimes
R(\tT)^{s_\alpha}\bigoplus K_\tT(\bar{T}^{+})\otimes
e^{\omega_{\alpha}(v)}R(\tT)^{s_{\alpha}}\ee as $K_\tT(\bar{T}^{+})\otimes
R(\tT)^W$-modules.

Note that, after localizing at $t_{\alpha}=\prod_{\beta\neq \alpha}
(1-e^{-{\beta(u)}})$, the only conditions defining $K_{\tG\times
  \tG}(X)$ in $\prod_{\sigma\in\cf_{+}(l)} R(\tT_{\sigma})\otimes
R(\tT)$ are the one corresponding to $\alpha$ and the condition
$(ii)$. This fact together with (\ref{ab2}), implies that
$K_{t_{\alpha}}\subseteq (M_{\alpha})_{t_{\alpha}}$.

Further by the following equalities \be\label{ab3}
R(\tT)^{s_{\alpha}}=\bigoplus_{\alpha\notin I} R(\tT)_{I}\ee
\be\label{ab4} e^{\omega_{\alpha}}\cdot
R(\tT)^{s_{\alpha}}=\bigoplus_{\alpha\in I} R(\tT)_{I}\ee and the
definition of $L$, it follows that $L_{t_{\alpha}}=
(M_{\alpha})_{t_{\alpha}}$. Moreover, since $L_{t_{\alpha}}\subseteq
K_{t_{\alpha}}\subseteq (M_{\alpha})_{t_{\alpha}}$, it follows that
$(K/L)_{t_{\alpha}}=0$ for every $\alpha\in\Delta$. \hfill $\Box$

Since the projective dimension of $(K/L)=1$, by Auslander-Buchsbaum
formula we know that $Supp~(K/L)$ is of pure codimension $1$ in
$\mbox{Spec} (R(\tT)\otimes R(\tG))$. Hence $Supp~(K/L)$ must contain a
prime ideal ${\mathfrak p}$ of height $1$ in $R(\tT)\otimes R(\tG)$. Since
$R(\tT)\otimes R(\tG)$ is a U.F.D, ${\mathfrak p}=(a)$ for some $a\in
R(\tT)\otimes R(\tG)$ and by Lemma \ref{loc} it follows that ${\mathfrak
  p}$ contains $1-e^{-{\alpha(u)}}$ and $1-e^{-{\beta(u)}}$ for
$\alpha\neq\beta\in\Delta$.

This implies that $a$ divides $1-e^{-{\alpha(u)}}$ and
$1-e^{-{\beta(u)}}$ for distinct $\alpha$ and $\beta$, which a
contradiction since $1-e^{-{\alpha(u)}}$ and $1-e^{-{\beta(u)}}$
are relatively prime in the U.F.D, $R(\tT)\otimes R(\tG)$ (see
\cite[ p.182]{bo}).

This contradiction implies that $K/L=0$ and hence $K=L$. This proves
(\ref{dsd}). We prove the remaining assertions of the theorem below.

Let $PL(\cf_{+})$ denote the piecewise linear functions on $\cc^{+}$
that are linear on the subdivision $\cf_{+}$. More precisely,
$PL(\cf_{+}):=$ \[\{h:\cc^{+}\ra \br~|~ \forall~ \sigma\in \cf_{+}(l)
~\mbox{and}~ v\in \sigma, ~h(v)=\langle h_{\sigma},v\rangle ~\mbox{for
  some} ~h_{\sigma}\in X^*(T)\}.\] By the parametrization of line
bundles on spherical varieties (see \cite[\S2.2]{br2}) we know that
the group of isomorphism classes of $\tG\times \tG$-linearised line
bundles on $X$ is isomorphic to $PL(\cf_{+})$. Let $\cl_{h}$ denote
the line bundle on $X$ corresponding to $h\in PL(\cf_{+})$.  Then
$\cl_{h}$ is $\tG\times \tG$-linearized and $\tB^{-}\times \tB$
acts on the fibre $\cl_{h}\mid_{z_{\sigma}}$ by the character
$(h_{\sigma},-h_{\sigma})$, where $z_{\sigma}$ denotes the base point
of the closed orbit $Z_{\sigma}$.  Moreover, $z_{\sigma}$ is also the
$T\times T$-fixed point in $\bar{T}^{+}$ corresponding to $\sigma\in
\cf_{+}(l)$.  Hence it follows that $L_{h}:=\cl_{h}\mid_{\bar{T}^{+}}$
is the $\tT\times \tT$-linearized line bundle on the toric variety
$\bar{T}^{+}$ corresponding to the piecewise linear function $h\in
PL(\cf_{+})$.  In particular, $\tT^{-}\times \tT$ acts on the fibre
$L_{h}\mid_{z_{\sigma}}$ by the character $h_{\sigma}$.

It suffices to show that the piece $K_{\tT}(\bar{T}^{+})\otimes 1$ in
(\ref{dsd}) is the subring of $K_{\tG\times \tG}(X)$ generated by
$\{[{\cl}_{h}]~:~h\in PL(\cf_{+})\}$, which is the subring generated
by  $Pic_{\tG\times \tG}(X)$.

Recall that under the canonical inclusion $K_{\tG\times \tG}(X)\hra
\prod_{\sigma\in\cf_{+}(l)} R({\tT}_{\sigma})\otimes R(\tT)$ obtained
via restriction to the base point $z_{\sigma}$ of the closed orbit
$Z_{\sigma}$, $\sigma\in \cf_{+}(l)$, $[\cl_{h}]$ maps to $(e^{h_{\sigma}}\otimes
1)_{\sigma\in \cf_{+}(1)}$.  Here we use the identification
$R(\tT\times \{1\})\otimes R(diag(\tT))\simeq R(\tT)\otimes R(\tT)$
coming from (\ref{ses}).

Also note that via restriction to $T$-fixed points $z_{\sigma}$,
$\sigma\in \cf_{+}(l)\}$, the image of $[L_h]\in
Pic^{\tT}(\bar{T}^{+})$ in $\prod_{\sigma\in
  \cf_{+}(l)}R({\tT}_{\sigma})$ is $(e^{h_{\sigma}})_{\sigma\in
  \cf_{+}(l)}$.  Moreover, it is also known that
$Pic^{\tT}(\bar{T}^{+})$ generates $K_{\tT}(\bar{T}^{+})$ as a ring
(see \cite[\S6.2]{vv}, where for $\rho\in \Delta(1)$, the piecewise
linear functions $u_{\rho}({\sigma})$ on the fan $\Delta$ which
generate $K_{T}(X(\Delta))$ as a ring, correspond to the classes of
the $T$-equivariant line bundles $L_{u_{\rho}}$). Thus it follows that
the image of $K_{\tT}(\bar{T}^{+})\otimes 1$ under restriction to
$\prod_{\sigma\in\cf_{+}(l)} R({\tT}_{\sigma})\otimes R(\tT)$ is same
as the image of the subring generated by $[\cl_{h}]$ for $h\in
PL(\cf_{+})$.  Hence the theorem.$\Box$

\bth\label{modkdec} We have the following isomorphism as
$K_{\tT}(\bar{T}^{+})\otimes R(\tT)^{W}$-submodules of
$\prod_{\sigma\in\cf_{+}(l)} R(\tT_{\sigma})\otimes R(\tT)$.
\[\bigoplus_{I\subseteq \Delta}\prod_{\alpha\in {I}}(1-e^{\alpha(u)})\cdot K_{\tT}(\bar{T}^{+})\otimes R(\tT)_{I}\simeq
\bigoplus_{I\subseteq \Delta}K_{\tT}(\bar{T}^{+})\otimes R(\tT)_{I}=
K_{\tG\times \tG}(X).\]

More explicitly, the above isomorphism maps an arbitrary element
$a\otimes b\in K_{\tT}(\bar{T}^{+})\otimes R(\tT)_{I}$ to the element
$\prod_{\alpha\in {I}}(1-e^{\alpha(u)})\cdot a\otimes b$. In particular, the
basis element $1\otimes f_v\in K_{\tT}(\bar{T}^{+})\otimes R(\tT)_{I}$ maps to
$\prod_{\alpha\in {I}}(1-e^{\alpha(u)})\otimes f_v\in K_{\tT}(\bar{T}^{+})\otimes R(\tT)_{I}$,
for $v\in C^{I}$ for every $I\subseteq \Delta$. 

Moreover, in $\bigoplus_{I\subseteq \Delta}K_{\tT}(\bar{T}^{+})\otimes
R(\tT)_{I}$ any two basis elements $1\otimes f_v$ and $1\otimes
f_{v^{\prime}}$ for $v\in C^{I}$, $v^{\prime}\in C^{I^{\prime}}$
($I,I^{\prime}\subseteq \Delta$) multiply as follows $(1\otimes
f_v)\cdot (1\otimes f_{v^{\prime}})=$
\be\label{modmult}\sum_{J\subseteq (I\cup I^{\prime})}\sum_{w\in
  C^{J}}(\prod_{\alpha\in I\cap I^{\prime}}(1-e^{\alpha(u)})\cdot
\prod_{\alpha\in (I\cup I^{\prime})\setminus
  J}(1-e^{\alpha(u)}))\otimes a^{w}_{v,v^{\prime}}\cdot (1\otimes
f_{w}).\ee \eeth

{\bf Proof:} Note that $\prod_{\alpha\in {I}}(1-e^{\alpha(u)})$ is not
a zero divisor in the integral domain $R(\tT)$ and hence in
$K_{\tT}(\bar{T}^{+})$ which is a free module over $R(\tT)$.  Thus each piece
$K_{\tT}(\bar{T}^{+}) \otimes R(\tT)_{I}$ is isomorphic to $\prod_{\alpha\in
  {I}}(1-e^{\alpha(u)})\cdot K_{\tT}(\bar{T}^{+})\otimes R(\tT)_{I}$, for every
$I\subseteq \Delta$, as $K_{\tT}(\bar{T}^{+})\otimes R(\tG)$-submodules of
$\prod_{\sigma\in\cf_{+}(l)}R({\tT}_{\sigma})\otimes R(\tT)$. The
isomorphism maps an element $a\otimes b\in K_{\tT}(\bar{T}^{+}) \otimes
R(\tT)_{I}$ to the element $\prod_{\alpha\in
  {I}}(1-e^{\alpha(u)})\cdot a\otimes b\in K_{\tT}(\bar{T}^{+})\otimes
R(\tT)_{I}$. This additively extends to an isomorphism of
$K_{\tT}(\bar{T}^{+})\otimes R(\tG)$-submodules of $K_{\tT}(\bar{T}^{+})\otimes R(\tT)$:
\be\label{iso}\bigoplus_{I\subseteq \Delta}\prod_{\alpha\in
  {I}}(1-e^{\alpha(u)})\cdot K_{\tT}(\bar{T}^{+})\otimes R(\tT)_{I}\simeq
\bigoplus_{I\subseteq \Delta}K_{\tT}(\bar{T}^{+})\otimes R(\tT)_{I}.\ee Further,
in the direct sum decomposition
$$K_{\tG\times \tG}(X)=\bigoplus_{I\subseteq
  \Delta} \prod_{\alpha\in {I}}(1-e^{\alpha(u)})\cdot K_{\tT}(\bar{T}^{+})\otimes
R(T)_{I},$$ the two basis elements $\prod_{\alpha\in
  {I}}(1-e^{\alpha(u)})\otimes f_{v}$ and $\prod_{\alpha\in
  {I'}}(1-e^{\alpha(u)})\cdot\otimes f_{v^{\prime}}$, $v\in C^{I}$
$v^{\prime}\in C^{I^{\prime}}$ multiply to give
\be\label{mbasis}\prod_{\alpha\in I\cap
  I^{\prime}}(1-e^{\alpha(u)})\cdot \prod_{\alpha\in I\cup
  I^{\prime}}(1-e^{\alpha(u)})\otimes f_v\cdot f_{v^{\prime}} .\ee

Note that the right hand side of (\ref{mbasis}) belongs to
\[\prod_{\alpha\in I\cup
  I^{\prime}}(1-e^{\alpha(u)})\cdot K_{\tT}(\bar{T}^{+})\otimes
R(T)^{W_{\Delta\setminus (I\cup I^{\prime})}}\subseteq
\bigoplus_{J\subseteq (I\cup I^{\prime})} \prod_{\alpha\in
  J}(1-e^{\alpha(u)})\cdot K_{\tT}(\bar{T}^{+})\otimes R(T)_{J}.\] Thus using
(\ref{mst}), (\ref{mbasis}) can be rewritten as 
\be\label{mbasis1}\sum_{J\subseteq (I\cup I^{\prime})}\sum_{w\in
  C^{J}}(\prod_{\alpha\in I\cap I^{\prime}}(1-e^{\alpha(u)})
\prod_{\alpha\in I\cup I^{\prime}}(1-e^{\alpha(u)})\otimes
a^{w}_{v,v^{\prime}})\cdot \prod_{\alpha\in J}(1-e^{\alpha(u)})\otimes
f_{w}).\ee

Now, using the isomorphism (\ref{iso}), in $\bigoplus_{I\subseteq
  \Delta}K_{\tT}(\bar{T}^{+})\otimes R(T)_{I}$ the basis elements $(1\otimes
f_v)$ and $(1\otimes f_{v'})$, $v\in C^{I}$ $v^{\prime}\in
C^{I^{\prime}}$, multiply as in (\ref{modmult}). \hfill $\Box$

\section{Ordinary K-ring of the regular compactification}

In this section we shall describe the ordinary $K$-ring of $X$ using
the results in Section 2.

Let $c_{K}:R(\tT)\ra K(G/B)$ denote the characteristic homomorphism
which sends $e^{\lambda}$ to $[\cl(\lambda)]$ for $\lambda\in
\Lambda$. In particular, $c_{K}(e^{\alpha})$ is $[\cl({\alpha})]$ for
$\alpha\in \Delta$. Further, we denote by $\bar{\lambda}_{I}\in
K(G/B)$, the image of $\prod_{\alpha\in I}(1-e^{-{\alpha}})\in R(\tT)$
for every $I\subseteq \Delta$.  Furthermore, let
$\bar{f_v}:=c_{K}(f_v)$ and $c^{w}_{v,v^{\prime}}\in \bz$ denote the
image under $c_{K}\mid_{R(\tG)}$ of $a^{w}_{v,v^{\prime}}\in R(\tG)$,
where $a^{w}_{v,v^{\prime}}$ is as in (\ref{mst}).

Further, if
\be \label{odec1} K(G/B)_{I}:=\bigoplus_{v\in C^{I}} \bz [{\bar f_v}]\ee then we have
\be\label{odec2} K(G/B)=\bigoplus_{I\subseteq \Delta}K(G/B)_{I}.\ee

Recall from \cite[Theorem 6.4]{vv} that we have the following
Stanley-Reisner presentation of the $\tT$-equivariant $K$-ring of $\bar{T}^{+}$:
\be\label{pres1} K_{\tT}(\bar{T}^{+})\simeq\frac{\bz[X_j^{\pm 1}:
  \rho_j\in\cf_{+}(1)]}{\langle X_F, ~~\forall~ F\notin \cf_{+}
\rangle}\ee where $X_{F}:=\prod_{\rho_i\in F}(1-X_i)$ for every
$F\subseteq \cf_{+}(1)$.

Further, since (\ref{rm}) is compatible with the canonical
$R(\tT)$-algebra structure on $K_{\tT}(\bar{T}^{+})$, and the $R(\tT)$-algebra
structure on $\prod_{\sigma\in\cf_{+}(l)}R(\tT_{\sigma})$ given by the
diagonal map, it follows that as an $R(\tT)$-algebra $K_{\tT}(\bar{T}^{+})$ has
the following presentation: \be\label{pres2}
K_{\tT}(\bar{T}^{+})\simeq\frac{R(\tT)[X_j^{\pm 1}: \rho_j\in\cf_{+}(1)]}{J}\ee
where $J$ is the ideal in $R(\tT)[X_j^{\pm 1}: \rho_j\in\cf_{+}(1)]$
generated by the elements $X_F$ for $F\notin \cf_{+}$ and
$(\prod_{\rho_j\in \cf_{+}(1)} X_{j}^{<\alpha,v_{j}>})-e^{\alpha}$ for
$\alpha\in \Delta$. Here $v_j$ denotes the primitive vector along the
edge $\rho_j$.

Let \be\label{deford}\cR(\bar{T}^{+}):=\bz\otimes_{R(\tG)}K_{\tT}(\bar{T}^{+})\ee where
$\bz$ is an $R(\tG)$-module under the augmentation
$\epsilon=c_{K}\mid_{R(G)}:R(\tG)\ra \bz$. Then it follows from
(\ref{pres2}) that \be\label{pres3ord} \cR(\bar{T}^{+})\simeq
\frac{K(G/B)[X_j^{\pm 1}: \rho_j\in\cf_{+}(1)]}{\mathfrak{J}}\ee where
$\mathfrak{J}$ is the ideal in $K(G/B)[X_j^{\pm 1}:
  \rho_j\in\cf_{+}(1)]$ generated by the elements $X_F$ for $F\notin
\cf_{+}$ and $(\prod_{\rho_j\in \cf_{+}(1)}
X_{j}^{<\alpha,v_{j}>})-[\cl({\alpha})]$ for $\alpha\in \Delta$. 

\bth\label{main}We have a canonical $\cR(\bar{T}^{+})$-module structure on
$K(X)$, induced from the $K_{\tT}(\bar{T}^{+})\otimes 1$-module structure on
$K_{\tG\times \tG}(X)$ given in Theorem \ref{modkdec}.  Moreover, $K(X)$
is a free module of rank $|W|$ over $\cR(\bar{T}^{+})$, $\cR(\bar{T}^{+})$ being
identified with the subring of $K(X)$ generated by $Pic(X)$.  More
explicitly, let \be\label{im1}\gamma_{v}:=1\otimes [{\bar
      f_v}]\in \cR(\bar{T}^{+})\otimes K(G/B)_{I}\ee for $v\in C^{I}$ for every
$I\subseteq \Delta$.  Then we have: \be\label{decord} K(X)\simeq
\bigoplus_{v\in W} \cR(\bar{T}^{+})\cdot \gamma_v .\ee Further, the above
isomorphism is a ring isomorphism, where the multiplication of any two
basis elements $\gamma_v$ and $\gamma_{v^{\prime}}$ is defined as
follows: \be\label{multord}\gamma_v\cdot
\gamma_{v^{\prime}}:=\sum_{J\subseteq (I\cup I^{\prime})}\sum_{w\in
  C^{J}}(\bar{\lambda}_{I\cap I^{\prime}}\cdot \bar{\lambda}_{(I\cup
  I^{\prime})\setminus J})\cdot c^{w}_{v,v^{\prime}} \cdot
\gamma_{w}.\ee \eeth

{\bf Proof:} By Theorem \ref{kdec} we have the following direct sum
decomposition of $K_{\tG\times \tG}(X)$ as an $R(\tT)\otimes
R(\tG)$-module \be\label{iso1}K_{\tG\times
  \tG}(X)\simeq\bigoplus_{I\subseteq
  \Delta}K_{\tT}(\bar{T}^{+})\otimes R(\tT)_{I}.\ee Now, using the
isomorphism \be\label{coe} K(X)\simeq K_{\tG\times
  \tG}(X)\otimes_{R(\tG)\otimes R(\tG)}\bz\ee (see \cite{mer}) and
(\ref{iso}) we get: \be\label{mainiso1} K(X)\simeq
\bigoplus_{I\subseteq \Delta} \cR(\bar{T}^{+})\otimes K(G/B)_{I}.\ee
Further, under the canonical restriction homomorphism \[K_{\tG\times
  \tG}(X)\ra K(X),\] the subring generated by $Pic^{\tG\times \tG}(X)$
in $K_{\tG\times \tG}(X)$, maps surjectively onto the subring
generated by $Pic(X)$ in $K(X)$. Hence by Theorem \ref{kdec} it
follows that in (\ref{mainiso1}), the subring generated by $Pic(X)$ in
$K(X)$ maps isomorphically onto the piece $\cR(\bar{T}^{+})\otimes
1\subseteq \cR(\bar{T}^{+})\otimes K_{\tT}(G/B)_{\emptyset}$.

Let $\gamma_{v}\in \cR(\bar{T}^{+})\otimes K_{\tT}(G/B)_I$ denote the element
$1\otimes \bar{f_v}$ for $v\in C^{I}$ for $I\subseteq \Delta$.  Then
by Theorem \ref{modkdec} and (\ref{mainiso1}), it follows that $K(X)$
is a free module of rank $|W|$ over the ring $\cR(\bar{T}^{+})$ with basis
$\gamma_v$ for $v\in W$, where $\cR(\bar{T}^{+})$ is identified with the subring
of $K(X)$ generated by $Pic(X)$.

 Further, recall from Theorem \ref{modkdec} that the multiplication of
 two basis elements $1\otimes f_{v}$ and $1\otimes f_{v^{\prime}}$ of
 $K_{\tG\times \tG}(X)\simeq\bigoplus_{I\subseteq \Delta}R(\tT)\otimes
 R(\tT)_{I}$ is defined as $(1\otimes f_v)\cdot (1\otimes
 f_{v^{\prime}})=$ \[\sum_{J\subseteq (I\cup I^{\prime})}\sum_{w\in
   C^{J}}\prod_{\alpha\in {I\cap I^{\prime}}} (1-e^{\alpha(u)})\cdot
 \prod_{\alpha\in (I\cup I^{\prime})\setminus
   J}(1-e^{\alpha(u)})\otimes a^{w}_{v,v^{\prime}}\cdot (1\otimes
 f_{w}).\] Thus their images $\gamma_v$ and $\gamma_{v^{\prime}}$ in $
 \bigoplus_{I\subseteq \Delta} \cR(\bar{T}^{+})\otimes K(G/B)_{I}$, multiply as
 follows:\be\label{mainmultn} \gamma_v\cdot
 \gamma_{v^{\prime}}:=\sum_{J\subseteq (I\cup I^{\prime})}\sum_{w\in
   C^{J}}(\lambda_{I\cap I^{\prime}})\cdot \lambda_{(I\cup
   I^{\prime})\setminus J})c^{w}_{v,v^{\prime}} \cdot \gamma_{w}.\ee

 We therefore conclude that the isomorphism \be\label{mainiso}
 K(X)\simeq \bigoplus_{v\in W} \cR(\bar{T}^{+}) \cdot \gamma_v \ee is in fact a
 ring isomorphism, where the multiplication of any two basis elements
 $\gamma_v$ and $\gamma_{v^{\prime}}$ for $v\in C^{I}$, $v^{\prime}\in
 C^{I^{\prime}}$ and $I,I^{\prime}\subseteq \Delta$ is defined as in
(\ref{mainmultn}).  \hfill$\Box$

\subsection{Relation with the equivariant $K$-ring of the wonderful
  compactification}

Recall from \cite[Proposition A2 ]{br} that we have a canonical
surjective morphism $f: X\rightarrow\bar{G_{ad}}$ from the regular
compactification of $G$ to the wonderful compactification of $G_{ad}$,
which further restricts to a morphism of toric varieties
$g:\bar{T}^{+}\rightarrow \bar{T_{ad}}^{+}$, which is also proper and
surjective. Furthermore, $f$ is equivariant with respect to the action
of $\tG\times \tG$, where $\tG\times \tG$ acts on $\bar{G_{ad}}$ via
the quotient $G^{ss}\times G^{ss}$.

In this section we relate using $f$, the $\tG\times \tG$-equivariant and
the ordinary $K$-ring of the regular compactification $X$, with the
equivariant and ordinary $K$-ring of the wonderful compactification
$\bar{G_{ad}}$ respectively.

It follows from (\ref{fc}) that $G^{ss}$ is the universal cover of
$G_{ad}$, and $T_{ad}:=T^{ss}/C$ is the maximal torus of
$G_{ad}$. Recall that ${rank}~(G_{ad})=rank~(G^{ss})=r$ which is the
semisimple rank of $G$.

Further, recall from \cite{dp} and \cite[\S2.2]{br2} that on
$\bar{G_{ad}}$ the isomorphism classes of line bundles correspond to
$\lambda\in X^*(T)$. Indeed, the line bundle $\cl_{\lambda}$ on $X$
associated to $\lambda$, admits a unique $G^{ss}\times
G^{ss}$-linearisation so that $(B^{ss})^{-}\times B^{ss}$ acts on the
fibre $\cl_{\lambda}\mid_{z}$ by the character $(\lambda,-\lambda)$,
where $z$ denotes the base point of the unique closed orbit
$G/B^{-}\times G/B$. In particular, $\cl_{\alpha_i}$ is $G^{ss}\times
G^{ss}$-linearised such that $(B^{ss})^{-}\times B^{ss}$ operates on
$\cl_{\alpha_i}\mid_{z}$, with the character $(\alpha_i,-\alpha_i)$, for
$1\leq i\leq r$. Moreover, $\cl_{\alpha_i}$ also admits a
$G^{ss}\times G^{ss}$-equivariant section $s_i$ whose zero locus is
the boundary divisor $D_i$ for $1\leq i\leq r$.

Further, since $\tG=G^{ss}\times C$, the line bundles $\cl_{\lambda}$
for $\lambda\in\Lambda$, in fact admit a $\tG\times
\tG$-linearization, by taking the $C\times C$ action to be trivial.

The ring $K_{\tG\times \tG}(X)$ gets the structure of an algebra over
the ring $K_{\tG\times \tG}(\bar{G_{ad}})$, by pull-back of
equivariant vector bundles along $f$.  The following theorem
explicitly describes this structure.

\bth\label{relwond} The ring $K_{\tG\times \tG}(X)$ can be described
as a $K_{\tG\times \tG}(\bar{G_{ad}})$-algebra as follows
\be\label{relwondeqn} K_{\tG\times \tG}(X)=\frac{K_{\tG\times
    \tG}(\bar{G_{ad}})[X_j^{\pm 1}:
  \rho_j\in\cf_{+}(1)]}{\mathfrak{J}}\ee where $\mathfrak{J}$ is the
ideal in $K_{\tG\times \tG}(\bar{G_{ad}})[X_j^{\pm 1}:
\rho_j\in\cf_{+}(1)]$ generated by the elements $X_F$ for $F\notin
\cf_{+}$ and $(\prod_{\rho_j\in \cf_{+}(1)}
X_{j}^{<\alpha,v_{j}>})-[\cl_{\alpha}]_{\tG\times \tG}$ for $\alpha\in
\Delta$. Here $[\cl_{\alpha}]_{\tG\times \tG}$ denotes the class of
the $\tG\times \tG$-equivariant line bundle $\cl_{\alpha}$ in
$K_{\tG\times \tG}(\bar{G_{ad}})$. \eeth

{\bf Proof:} Note that $\tG\times \tG$ acts on $\bar{G_{ad}}$ via the
quotient $G^{ss}\times G^{ss}$, so that $\tC\times \tC$ acts trivially
on $\bar{G_{ad}}$. We therefore have the following isomorphism of
$R(\tG)\otimes R(\tG)$-algebras \be\label{e4} K_{\tG\times
  \tG}(\bar{G_{ad}})=K_{G^{ss}\times G^{ss}}(\bar{G_{ad}})\otimes
R(\tC\times\tC).\ee From \cite[Theorem 3.3]{u}, using the equalities
(\ref{e1}), (\ref{e2}), (\ref{e3}), (\ref{2}) and (\ref{e4}), it
follows that the $\tG\times \tG$-equivariant $K$-theory of
$\bar{G_{ad}}$ has the following description as an $R(\tT)\otimes
R(\tG)$-module \be\label{wd} K_{\tG\times
  \tG}(\bar{G_{ad}})=\bigoplus_{I\subseteq \Delta}\prod_{\alpha\in
  I}(1-e^{-{\alpha(u)}}) R(\tT)\otimes R(\tT)_{I}.\ee (Note that in
\cite[Section 3]{u}, $G$ and $T$ denote $G^{ss}$ and $T^{ss}$
respectively and $X$ denotes $\bar{G_{ad}}$.) We now see that
(\ref{dsd}) and (\ref{wd}) together imply the following isomorphism of
$K_{\tT}(\bar{T}^{+})\otimes 1$-algebras \be\label{forregwond}
K_{\tG\times \tG}(X)=K_{\tG\times
  \tG}(\bar{G_{ad}})\otimes_{R(\tT)\otimes
  1}(K_{\tT}(\bar{T}^{+})\otimes 1),\ee where the $R(\tT)\otimes
1$-algebra structure on $K_{\tG\times \tG}(\bar{G_{ad}})$ is from the
first factor in the direct sum decomposition (\ref{wd}). Indeed,
$R(\tT)\otimes 1$ is identified with the subring of $K_{\tG\times
  \tG}(\bar{G_{ad}})$ generated by $Pic_{\tG\times
  \tG}(\bar{G_{ad}})$. Under this identification, the isomorphism
class of $\cl_{\alpha}$ in $K_{\tG\times \tG}(\bar{G_{ad}})$
corresponds to $e^{\alpha}\otimes 1$ in $R(\tT)\otimes 1$ for
$\alpha\in \Delta$ (see \cite[Theorem 3.6]{u}).  Moreover, since
$K_{\tT}(\bar{T}^{+})$ has the presentation given by (\ref{pres2}) as
an $R(\tT)$-algebra, the theorem follows.  $\Box$

\subsubsection{Relation with the ordinary $K$-ring of the wonderful
  compactification}

The ring $K(X)$ gets the structure of an algebra over the ring
$K(\bar{G_{ad}})$, via pull-back of vector bundles along $f$.  The
following theorem explicitly describes this structure.

\bth\label{relwond1} The ring $K(X)$ has the following presentation as
a $K(\bar{G_{ad}})$-algebra: \be\label{forregwond1}
K(X)=\frac{K(\bar{G_{ad}}) [X_j^{\pm 1}:
  \rho_j\in\cf_{+}(1)]}{\mathfrak{I}}\ee where $\mathfrak{I}$ is the
ideal in $K(\bar{G_{ad}})[X_j^{\pm 1}: \rho_j\in\cf_{+}(1)]$ generated
by the elements $X_F$ for $F\notin \cf_{+}$ and $(\prod_{\rho_j\in
  \cf_{+}(1)} X_{j}^{<\alpha,v_{j}>})-[\cl_{\alpha}]$ for $\alpha\in
\Delta$. Here $[\cl_{\alpha}]$ denotes the class of the line bundle
$\cl_{\alpha}$ in $K(\bar{G_{ad}})$.  \eeth

{\bf Proof:} Using the fact that $\tC\times \tC$ acts trivially on
$\bar{G_{ad}}$, by (\ref{eqk3}) we have \be\label{owd1}
K_{G^{ss}\times G^{ss}}(\bar{G_{ad}})\otimes_{R(G^{ss}\times
  G^{ss})}\bz=K_{\tG\times \tG}(\bar{G_{ad}}) \otimes_{R(\tG\times
  \tG)}\bz=K(\bar{G_{ad}}).\ee Again by (\ref{eqk3}), \be\label{oreg}
K_{\tG\times \tG}(X)\otimes_{R(\tG\times \tG)}\bz=K(X).\ee Moreover,
note that $[\cl_{\alpha}]_{\tG\times \tG}$ restricts to
$[\cl_{\alpha}]$ under the canonical forgetful homomorphism from
$Pic_{\tG\times \tG}(\bar{G_{ad}})$ to $Pic(\bar{G_{ad}})$. The
equality (\ref{forregwond1}) now follows from (\ref{relwondeqn}),
(\ref{owd1}) and (\ref{oreg}).$\Box$

Consider the principal $T_{ad}$-bundle $P$ corresponding to the direct
sum of the line bundles $\cv:=\bigoplus_{1\leq i\leq r}
\cl_{\alpha_{i}}$ on $\bar{G_{ad}}$ (see \cite[Section 5]{dp1} and
\cite[Section 10]{bdp}). Then $P$ can be thought of as a principal
$T$-bundle where $T$ acts via its quotient $T_{ad}$. Then
$P\times_{T}\bar{T}^{+}$ is a toric bundle with fibre $\bar{T}^{+}$.
Moreover, since each $\cl_{\alpha_i}$ is $\tG\times \tG$-linearized,
$P$ is a left $\tG\times \tG$-space where the bundle map $\pi:P\ra
\bar{G_{ad}}$ is $\tG\times \tG$-equivariant for the canonical
$\tG\times \tG$-action on $\bar{G_{ad}}$. Moreover, the right
$T$-action is compatible with the left $\tG\times \tG$-action on
$P$.

\bcor\label{eqtorbund} The ring $K_{\tG\times\tG}(X)$ as a
$K_{\tG\times \tG}(\bar{G_{ad}})$-algebra and $K(X)$ as a
$K(\bar{G_{ad}})$-algebra are respectively isomorphic to the
$\tG\times \tG$-equivariant and ordinary Grothendieck ring of the
above defined toric bundle $P\times_{T} \bar{T}^{+}$ over
$\bar{G_{ad}}$. \ecor

{\bf Proof:} Since $X$ is a projective regular embedding of $G$ and $\bar{T}^{+}$
is the inverse image of $\ba^r$ under the canonical morphism $f:X\ra
\bar{G_{ad}}$, the restriction $g:\bar{T}^{+}\ra \ba^r$ of the
projective morphism $f$ is a projective morphism of toric
varieties. This implies in particular that $\bar{T}^{+}$ is a
semi-projective $T$-toric variety.

Hence by Theorem \ref{ektb} of \S4 it follows that the right hand side
of (\ref{relwondeqn}) and (\ref{forregwond1}) in Theorem \ref{relwond}
and Theorem \ref{relwond1} are isomorphic respectively to the
$\tG\times \tG$-equivariant and ordinary Grothendieck ring of the
toric bundle $P\times_{T} \bar{T}^{+}$ over $\bar{G_{ad}}$. $\Box$

\brem Using techniques similar to \cite[Section 4]{u1} we can
alternately describe $K(X)_{\bq}$ as \be\label{altiso}
K(X)_{\bq}\simeq \bigoplus_{v\in W} \cR(\bar{T}^{+}) \cdot \gamma_v,
\ee where $\gamma_v=[\co_{X^v}]\in K(G/B)_{\bq}$ for $v\in
C^{I}$. Moreover, the multiplication of $\gamma_v$ and
$\gamma_{v^{\prime}}$ for $v\in C^{I}$, $v^{\prime}\in C^{I^{\prime}}$
and $I,I^{\prime}\subseteq \Delta$ is as defined in (\ref{mainmultn}),
where $c^w_{v,v'}$ are now the multiplicative structure constants of
the opposite Schubert classes described for instance in
\cite[(3.82)]{u1}. Since this description involves the Schubert classes,
it is more geometric in nature.  However, we require rational
coefficients for this description, since we use the lifts of the
Schubert classes in $R(\tT)$ in the place of the Steinberg basis, and
these lifts form a basis for $R(\tT)$ over $R(\tG)$ only after
localizing at the augmentation ideal $I(\tG)=\{a-\epsilon(a):a\in
R(\tG)\}$.  \erem

\section{Appendix}

\subsection{Semi-projective toric varieties}

In this section we recall the geometric definition and some essential
properties of semi-projective toric varieties. Our description here is
brief mainly fixing the notations and conventions necessary to state
and prove the main result Theorem \ref{ektb} in the next section (for
details we refer to \cite[Section 2]{hs}).

\subsubsection{The geometric definition}
For every scheme $X$ there is a canonical morphism $\pi_{X}: X\ra X_0$
to the affine scheme $X_0:=\mbox{Spec}(H^0(X,\co_X))$ of regular
functions on $X$.  We call a toric variety $X$ {\it semi-projective}
if $X$ has at least one torus-fixed point and the morphism $\pi_X$ is
projective (see \cite[p. 501]{hs}).

There are equivalent characterizations of a semi-projective toric
variety in terms of the combinatorics of its defining fan, and also
via geometric invariant theory (see \cite[Theorem 2.6 and Corollary
2.7]{hs}).

\subsubsection{Combinatorial characterization}

In this section we describe the fan associated to a smooth
semi-projective toric variety \cite[p. 499]{hs}.

Let $\cb=\{v_1,\ldots,v_d \}$ be a configuration of vectors in the
lattice $N\simeq \bz^n$. Further, let $pos(\cb)$ be the convex
polyhedral cone in $N_{\br}:=N\otimes_{\bz}\br$ spanned by the vectors
in $\cb$. A {\it triangulation} of $\cb$ is a simplicial fan $\Delta$
in $N$, whose rays lie in $\cb$ and whose support equals
$pos(\cb)$. Moreover, the triangulation is said to be {\it unimodular}
if every maximal cone of $\Delta$ is spanned by a basis of $N$. This
is equivalent to the toric variety $X(\Delta)$ being smooth.

A {\it $T$-Cartier divisor} on a fan $\Delta$ is a continuous function
$\Psi:pos(\cb)\ra \br$ which is linear on each cone of $\Delta$ and
takes integer values on $N\cap pos(\cb)$. It is further said to be
{\it ample} if the function $\Psi:pos(\cb)\ra \br$ is convex and
restricts to a different linear function on each maximal cone of
$\Delta$. Indeed, $\Psi$ is the piecewise linear support function
corresponding to an ample line bundle on $X(\Delta)$.

We say that the {\it triangulation} $\Delta$ of $\cb$ is {\it
  regular} if there exists a $T$-Cartier divisor on $\Delta$ which is
{\it ample}.

By \cite[Corollary 2.7]{hs} which characterizes smooth semi-projective
toric varieties, we recall that its associated fan $\Delta$ is a {\it
  regular unimodular triangulation} of a configuration of vectors
$\cb$ which spans the lattice $N$.

\subsubsection{Characterization via GIT}

Let $\Delta$ be a regular unimodular triangulation of $\cb$ as above
so that $X(\Delta)$ is a smooth semi-projective toric variety.  

In this section we first outline the construction by which we realize
$X(\Delta)$ as a projective GIT quotient of the complex affine space
$\bc^d$ by $(\bc^*)^{d-n}$, which is an $(d-n)$-dimensional subtorus
of $(\bc^*)^d$. (see \cite[Theorem 2.4]{hs}).

Let $M:={\cx}^{*}(T)=Hom(N,\bz)$ be the dual lattice of characters of
$T$ and $M_{\br}:=M\otimes \br$. We have the $d\times n$ integer
matrix \be\label{mat} \mathcal{Q}:=[v_1,\ldots,v_d]^T.\ee

Consider the map $M\ra \bz^d$ induced by the matrix
$\mathcal{Q}$. Since $\Delta$ is smooth and has at least one maximal
dimensional cone it follows that the cokernel of this map is
identified with $\bz^{d-n}$ and is isomorphic to the Picard group of
$X(\Delta)$ (see \cite[Proposition on p. 63]{f}). We therefore have an
exact sequence \be\label{es1}0\ra M\ra \bz^d\ra Pic(X) \ra 0.\ee
Applying $\mbox{Hom}(.,\bc^*)$ to (\ref{es1}) gives us the exact
sequence \be\label{es2}1\ra G:=(\bc^*)^{d-n}\ra (\bc^*)^d\ra
(\bc^*)^n\ra 1.\ee Thus the group $G$ is embedded as a
$(d-n)$-dimensional subtorus of $(\bc^*)^d$ and hence acts on the
affine space $\bc^d$. Further, let $\theta:=[L_{\Psi}]$ denote the
class in $Pic(X)$ of the ample line bundle $L_{\Psi}$ corresponding to
the support function $\Psi$.

We assume the action of $G$ on $\bc^d$ to be linearized on the trivial
bundle via the ample character $\theta$ of $G$. Then the geometric
quotient of the semi-stable locus of $\bc^d$ modulo $G$ coincides with
the toric variety $X(\Delta)$ (see \cite[Section 2]{hs} or \cite[Section 7.2 and
Section 14.2]{cls} for details).

Conversely any projective GIT quotient of the affine space $\bc^d$
modulo a subtorus $G$ of $(\bc^*)^d$, with respect to a
$G$-linearization on the trivial bundle on $\bc^d$ corresponding to a
character $\theta$ of $G$, gives a semi-projective toric variety.

Indeed applying $Hom(.,\bc^*)$ to the inclusion $G \hra (\bc^*)^d$ and
using the fact that $\bc^*$ is divisible, we get the following
surjective map of abelian groups $\bz^d\twoheadrightarrow
\bz^{d-n}=Hom(G,\bc^*)$. Let $a_i$ be the image in $\bz^{d-n}$ of the
standard basis element $e_i$ of $\bz^d$ under this surjection. We then
have the vector configuration $\ca:=\{ a_1,\ldots, a_d\}$ in
$\bz^{d-n}$.

Let $S=\bc[x_1,x_2,\ldots,x_d]$ be the ring of polynomial functions on
$\bc^d$. Then $S$ admits a grading by the semigroup $\bn\ca\subseteq
\bz^{d-n}$, setting $deg(x_i)=a_i$. Further, we have a natural $G$
action on $S$ defined by $g\cdot x_i:=a_i(g)\cdot x_i$ for $g\in G$,
$1\leq i\leq d$.  A polynomial in $S$ is homogeneous if and only if it
is a $G$-eigenvector. Let $S_{\theta}$ denote the homogeneous
polynomials of degree $\theta$. Then $S_{\theta}$ is a module over the
subalgebra $S_0$ of degree $0$ polynomials. Further, the ring
$S_{(\theta)}:=\oplus_{r=0}^{\infty}S_{r\theta}$ is a finitely
generated $S_0$-algebra. The projective GIT quotient
$\bc^d//_{\theta}G=\mbox{Proj} (S_{(\theta)})$ is then a toric variety
that is projective over the affine toric variety
$\bc^d//_0G=\mbox{Spec}(S_0)$ (see \cite[Section 14.2]{cls} and
\cite[Definition 2.2]{hs}).

\subsection{Equivariant $K$-theory of semi-projective toric varieties}

In this section we shall show that the result of \cite[Theorem
1.2(iv)]{su} on the structure of the Grothendieck ring of a toric
bundle associated to smooth projective toric varieties, holds for
toric bundles whose fibre is a smooth semi-projective toric variety.

We consider the following setting which is more general than the one
considered in \cite{su}.

Let $\pi: E\rightarrow B$ be an algebraic principal $T$-bundle, where
$B$ is an irreducible, non-singular and noetherian scheme over
$\bc$. Let $G$ be a connected complex reductive algebraic group. We
let $E$ and $B$ be algebraic left $G$-spaces with $\pi$ being
$G$-invariant, so that the left $G$-action on $E$ commutes with the
right $T\simeq (\bc^{*})^l$-action.

Let $X$ be a smooth semi-projective $T$-toric variety and $\Delta$ be
its associated fan in $N$. Without loss of generality, for simplicity
of notation we assume that $\{v_1,\ldots,v_d\}$ is the set of
primitive vectors along the set of edges
$\Delta(1):=\{\rho_1,\ldots,\rho_d\}$ of $\Delta$.

Let $V(\gamma)$ denote the orbit closure of the $T$-orbit $O_{\gamma}$
of the toric variety $X$ corresponding to a cone
$\gamma\in\Delta$. Let $U_{\sigma}$ denote the $T$-stable open affine
subvariety corresponding to a cone $\sigma\in \Delta$.

Further, corresponding to every edge $\rho_j$, we have a canonical
$T$-equivariant line bundle $L_j$ on $X$. Moreover, on $X$ there are
trivial line bundles $L_u=X\times \bc_{u}$ corresponding to characters
$u\in X^*(T)$. Here the bundle projection $L_u\stackrel{h_u}{\ra} X$
is $T$-equivariant where $T$ acts on the fibre $\bc_{u}$ via the
character $\chi_u:T\ra \bc^*$.  (see \cite[Chapter 3]{f}).

Consider the associated toric bundle $E(X):= E\times_{T}
X$ where $X:=X(\Delta)$ is a smooth semi-projective $T$-toric
variety associated to the fan $\Delta$.

We define $E(L_j):=E\times_{T}L_j$, the line bundle on $E(X)$
associated to $L_j$ for $1\leq j\leq d$. Also, let
$E(L_u)=E\times_{T}L_u$ be the line bundle on $E(X)$ associated to
$L_u$ with the bundle map $[e,v]\mapsto [e,h_{u}(v)] $. Further, both
$E(L_j)$ and $E(L_u)$ are both line bundles on $B$ with the bundle map
$[e,v]\mapsto \pi(e)$. Furthermore, $E(L_j)$ as well as $E(L_u)$ are both
$G$-linearized, where the $G$-linearization comes from the left
$G$-action on $E$.

Associated to the $T$-stable subvarieties $V(\gamma)$ corresponding to
a cone $\gamma$ in $\Delta$, we define the subvariety
$E(V(\gamma)):=E\times_T V(\gamma)$ of $E(X)$.

Let $x$ be a $T$-fixed point in $X$. Then the projection $p:E(X)\ra B$
has a canonical section $s:B\ra E(X)$ defined by $s(b)=[e,x]$ where
$e\in \pi^{-1}(b)$, so that $p\circ s=id_{B}$. Thus the corresponding
induced maps $s^*:K_{G}(E(X))\ra K_{G}(B)$ and $p^*:K_{G}(B)\ra
K_{G}(E(X))$ of the $G$-equivariant Grothendieck rings satisfy
$s^*\circ p^*=id^*_{B}$. In particular, $p^{*}$ is injective and
$K_{G}(E(X))$ gets the structure of $K_{G}(B)$-module via $p^{*}$.

We have the following theorem on the $G$-equivariant Grothendieck ring
of the toric bundle $E(X)$.

\bth\label{ektb} The $G$-equivariant $K$-ring of $E(X)$ has the
following presentation as a $K_{G}(B)$-algebra: \be\label{srpres}
K_{G}(E(X))=\frac{K_{G}(B)[X_{j}^{\pm1}]}{\langle X_{F}~:~
  F\notin\Delta~~;~~\prod_{\rho_j\in\Delta(1)} (1-X_{j})^{<u,v_j>}
  -[E(L_{u})]_{G},~\forall~u\in M\rangle}\ee where
$X_{F}=\prod_{\rho_j\in F}(1-X_j)$. (Since $E(L_{u})\ra B$ is a
$G$-linearized line bundle on $B$ (as observed above), we can
consider its class $[E(L_u)]_{G}\in K_{G}(B)$. )\eeth

We now fix some more notations and also recall necessary details about
the structure of $X$ and $\Delta$ before we go to the proof of the
above theorem.

The $T$-fixed locus in $X$ consist of the set of $T$-fixed
points \be\label{fp}\{x_{\sigma_1},
x_{\sigma_2}\ldots,x_{\sigma_m}\}\ee corresponding to the set of
maximal dimensional cones
\be\label{maxc}\Delta(n):=\{\sigma_1,\sigma_2,\ldots,\sigma_m\}.\ee

Choose a {\it generic} one parameter subgroup $\lambda_v\in X_{*}(T)$
corresponding to a $v\in N$ which is outside all hyperplanes spanned
by $(n-1)$-dimensional cones, so that the fixed points of $\lambda_v$
is the set $\{x_{\sigma_1},\ldots,x_{\sigma_m}\}$ (see
\cite[\S3.1]{br3}).

For each $x_{\sigma_i}$, we define the plus strata as
follows \be\label{ps} Y_i=\{x\in X \mid \lim_{ t\ra 0} \lambda_{v}
(t)x~ \mbox{exists and is equal to}~ x_{\sigma_i}\}.\ee

Then \be\label{pu} Y_i=\bigcup_{\tau_i\subseteq \gamma\subseteq
  \sigma_i}O_{\gamma}.\ee The union here is over the set of faces
$\gamma$ of $\sigma$ satisfying the property that the image of $v$ in
$N_{\br}/\br\gamma$ is in the relative interior of
$\sigma_i/\br\gamma$ (see \cite[Lemma 2.10]{hs}). Since the set of
such faces is closed under intersections, we can choose a minimal such
face of $\sigma_i$ which we denote by $\tau_i$. 

Further, if we choose $v\in|\Delta|$, then we have $X=\bigcup_{i=1}^m
Y_i$, which is the Bialynicki-Birula decomposition of the toric
variety $X$ with respect to the one-parameter subgroup $\lambda_v$.
From (\ref{pu}) it follows that \be\label{paff}Y_i=V(\tau_i) \cap
U_{\sigma_i}\simeq \bc^{n-k_i},\ee where $k_i=\dim(\tau_i)$ for $1\leq
i\leq m$. Indeed, $Y_i$ is an affine open set in the toric variety
$V(\tau_i)$ corresponding to a maximal $(n-k_i)$-dimensional cone in
$N/\br \tau_i$ (see \cite[Chapter 3]{f}).

Corresponding to the ample $T$-Cartier divisor $\Psi$ on $\Delta$, we
have the moment map $\mu:X\ra \br^n$ (see \cite[p.503]{hs}). Further,
the image of $X$ under $\mu$ is identified with a convex
polyhedron \[Q_{\Psi}:=\{v\in M_{\br}:{\mathcal{Q}}v\geq \psi\}\] (see
\cite[p.500]{hs}) where $\psi:=(\Psi(v_1),\ldots,\Psi(v_d))\in \bz^d$
and $\mathcal Q$ as in (\ref{mat}).

Also, we can define the moment map $\mu_v$ for the circle action on
$X$ induced by the {\it generic} one-parameter subgroup $\lambda_v$ by
$\mu_v(x):=\langle v,\mu(x)\rangle$. Further, using $\mu_v$, we can
define an ordering \be\label{o1}\sigma_1,\sigma_2,\ldots,\sigma_m\ee
of the maximal dimensional cones so that\be\label{mo}
\mu_v(x_{\sigma_i})\leq \mu_v(x_{\sigma_j})~~
\mbox{implies}~~\sigma_i\leq \sigma_j.\ee

Then the distinguished faces $\tau_i\subseteq \sigma_i$, with respect
to this ordering satisfy the following $(*)$ property: \be\label{star}
\tau_i\subseteq \sigma_j~~\mbox{implies}~~i\leq j.\ee This is
equivalent to \be\label{clstar}\bar{Y_{i}}\subseteq \bigcup_{j\geq i}
Y_{j}\ee for every $1\leq i\leq m$. This means that the
Bialynicki-Birula decomposition of $X$ is filtrable in the sense of
\cite[\S3.2]{br3}. We can order the cells $Y_1, Y_2,\ldots, Y_m$ so
that (\ref{clstar}) holds or equivalently, each $Y_i$ is closed in
$Y_1\cup Y_2\cup \cdots \cup Y_{i}$.

Let \[Z_i:=\bigcup_{j\geq i} Y_j.\] By (\ref{clstar}) we see that
$Z_i$ are closed subvarieties, which form a chain
\be\label{acd}X=Z_1\supseteq Z_2\supseteq \cdots \supseteq
Z_m=\{x_{\sigma_m}\}, \ee such that $Z_{i}\setminus Z_{i+1}=Y_i\simeq
\bc^{n-k_i}$. This defines an algebraic cell decomposition on $X$. In
other words $X$ is paved by the affine spaces $Y_i$ for $1\leq i\leq
m$.

Let $V_i:=V(\tau_i)$ corresponding to the distinguished faces $\tau_i$
of the maximal cone $\sigma_i$ for $1\leq i\leq m$. Let $E(V_i):= E\times_{T} V_i$ for
$1\leq i\leq m$. Also $E(Y_i):=E\times_{T} Y_i$ and
$E(Z_i):=E\times_{T} Z_i$.

We now prove the main theorem:

{\bf Proof:} Let $\cR=\cR(K_{G}(B),\Delta)$ denote the ring defined in
\cite[Definition 1.1]{su}, where $r_i=1- [E(L_{u_i})]_{G}$ in
$K_{G}(B)$ corresponding to the basis of characters $u_i\in X^{*}(T)$,
$1\leq i\leq n$ dual to a basis $v_1,\ldots, v_n$ of $N$ (we can
renumber the edges so that the first $n$ primitive edge vectors span a
maximal cone and hence form a basis of $N$).

We have a canonical homomorphism $\psi:\cR\ra K_{G}(E(X))$ of
$K_{G}(B)$ algebras defined by $x_j\ra 1-[E(L_j)]$ for every $1\leq
j\leq d$. Now, by \cite[Lemma 2.2 (iv)]{su} we know that $x(\tau_i)$,
$1\leq i\leq m$ span $\cR$ as a $K_{G}(B)$-module.  Thus if we show
that $K_{G}(E(X))$ is a free $K_{G}(B)$-module of rank $m$ with basis
$[\co_{E(V_i)}]$, $1\leq i\leq m$ (which are respectively the images
of $x(\tau_i)$, $1\leq i\leq m$ under $\psi$) then it will follow that
$\psi$ is an isomorphism and we are done.  

{\it Claim:} The classes $[\co_{E(V_i)}]$ $1\leq i\leq m$, form a basis
for $K_{G}(E(X))$ as a $K_{G}(B)$-module.

{\it Proof of Claim:} Firstly, from \cite[Corollary 6.10(i)]{vv} and
(\ref{pres2}) above it follows that $K(X)$ is isomorphic to the ring
$\cR(\bz,\Delta)$, with $r_i=1$ for $1\leq i\leq n$. Moreover,
$x(\tau_i)$ maps to $[\co_{V_i}]$ for $1\leq i\leq m$ under this
isomorphism.  Now, by \cite[Lemma 2.2 (iv)]{su} we know that
$x(\tau_i)$, $1\leq i\leq m$ spans $\cR(\bz,\Delta)$ as a
$\bz$-module. Thus $K(X)$ is spanned as a $\bz$-module by
$[\co_{V_i}]$, $1\leq i\leq m$. However, by \cite[Corollary
6.10(ii)]{vv} we also know that $K(X)$ is a free abelian group of rank
$m$, which is the number of maximal cones in $\Delta$ (this also
follows by cellular fibration lemma \cite[Lemma 5.5.1]{cg} since $m$
is the number of cells in the algebraic cell decomposition (\ref{acd})
of $X$). Hence $[\co_{V_i}]$, $1\leq i\leq m$ form a basis for $K(X)$
as a free $\bz$-module.

We now define $\phi:K_{G}(B)\otimes K(X)\ra K_{G}(E(X))$ by
\be\sum_{1\leq i \leq m}b_i\otimes [\co_{V_i}]\mapsto\sum_{1\leq i\leq
  m}p^{*}(b_i)[\co_{E(V_i)}].\ee Now, it remains only to show that the
$K_{G}(B)$-module homomorphism $\phi$ is injective and
surjective. This we prove below.

Note that we have the filtration \be\label{assald}E(Z_1)=E(X)\supseteq
E(Z_2)\supseteq\cdots \supseteq E(Z_m)\simeq B\ee where
$E(Z_i)\setminus E(Z_{i+1})=E(Y_i)$ is an affine bundle on
$B$. Indeed, the restriction of $p:E(X)\ra B$ to $E(Z_i)$ and $E(Y_i)$
are also $G$-equivariant with respect to the left $G$-action on $E$
and $B$. Then by the cellular fibration lemma \cite[Lemma 5.5.1]{cg},
for every $1\leq i\leq m$ we have a short exact sequence
\be\label{cfses}0\ra K_{G}(E(Z_{i+1}))\ra K_{G}(E(Z_i))\ra
K_{G}(E(Y_i))\ra 0.\ee We claim that the restriction of $\phi$ defines
isomorphisms $\phi_i:K_{G}(B)\otimes K(Z_i)\ra K_{G}(E(Z_i))$ for each
$1\leq i\leq m$. We prove this by downward induction on $i$. This is
trivially true for $i=m$, since in this case $E(Z_m)\simeq
B$. Consider the commutative diagram
\be\label{diagramcf} \begin{array}{llllllllll}& & K_{G}(B)\otimes
  K(Z_{i+1}) &\ra& K_{G}(B)\otimes K(Z_i)&\ra& K_{G}(B)\otimes
  K(Y_i)&\ra 0\\ & &\hspace{1.4cm}\da\phi_{i+1} & &\hspace{1.4cm}\da\phi_{i} &
  & \hspace{1.4cm}\da & & \\ & 0\ra
  &\hspace{0.7cm}K_{G}(E(Z_{i+1}))&\ra&\hspace{0.7cm}K_{G}(E(Z_i))&\ra&
  \hspace{0.7cm}K_{G}(E(Y_i)) \end{array}\ee where the bottom row is a part of
(\ref{cfses}) and the top horizontal row is obtained from taking
tensor product of the exact sequence \be\label{cfses1}0\ra
K(Z_{i+1})\ra K(Z_i)\ra K(Y_i)\ra 0\ee with $K_{G}(B)$.  The exact
sequence (\ref{cfses1}) is a consequence of the cellular fibration
lemma in the case when the base $B$ is a point. Moreover, by the Thom
isomorphism \cite[Theorem 5.4.11]{cg} $p^*: K_{G}(B)\simeq
K_{G}(E(Y_i))$ is an isomorphism. Furthermore, since $K(Y_i)\simeq
\bz$ (see for example \cite{gu}), it follows that the homomorphism
$K_{G}(B)\otimes K(Y_i)\ra K_{G}(E(Y_i))$ in (\ref{diagramcf}) is an
isomorphism. Therefore if we assume the $\phi_{i+1}$ is surjective
then it follows by diagram chase that $\phi_i$ is
surjective. Similarly if we assume that $\phi_{i+1}$ is injective,
then since the bottom horizontal row is left exact by (\ref{cfses}),
it follows again by diagram chase that $\phi_i$ is injective.  $\Box$


\begin{thebibliography}{9}
\bibitem{bdp} E. Bifet, C. De Concini and C. Procesi, {\it Cohomology of
regular embeddings}, Adv. Math. {\bf 82} (1990), 1-34.


\bibitem{bo} N.Bourbaki, {\it Groupes et alg\`{e}bres de Lie,}
Chapitres 4, 5 et 6.- Paris, Hermann, 1968 ({\it Act. scient. et ind.,
1337; Bourbaki,} 34).


\bibitem{br} M. Brion, {\it The behaviour at infinity of the Bruhat
decomposition},  Comment. Math. Helv. {\bf 73} (1998), 137-174.


\bibitem{br2} M. Brion, {\it Groupe de Picard et nombres caracteristiques
des varietes spheriques}, {Duke Math. J.}{\bf 58}, No. 2, 397-424,
1989.

\bibitem{br3} M. Brion, {\it Equivariant Chow groups for torus actions},
  Transformation Groups {\bf 2}, No. 3, 225-267, 1997.

\bibitem{cg} N.Chriss and V.Ginzburg, {\it Representation theory and 
complex geometry}, Birkhauser (1997).

\bibitem{cls}  D. A. Cox, J. B. Little and H.K.Schenck, {\it Toric
    Varieties},  Graduate Studies in Mathematics {\bf 124}, American
  Mathematical Society, 2011. 


\bibitem{dp} C. De Concini and C. Procesi, {\it Complete symmetric
varieties}, Invariant Theory (Proceedings, Montecatini 1982),
Lecture Note in Math. {\bf 996}, 1-44, Springer-Verlag, New York
1983.

\bibitem{dp1}C. De Concini and C. Procesi, {\it Complete symmetric
  varieties, II: Intersection theory},  Algebraic Groups and related
  topics (R. Hotta, ed.), Adv. Studies in Pure Math., {\bf 6}, pp
  481-513, North Holland, Amsterdam, 1985. 


\bibitem{f} W. Fulton, {\it Introduction to toric varieties}, Annals
  of Mathematical Studies {\bf 131}, Princeton University Press,
  Princeton, New Jersey (1993). 


\bibitem{gu} J. Gubeladze, {\it K-theory of affine toric varieties}, 
Homology, Homotopy and Applications {\bf 1} (1999), 135-145.




\bibitem{h} J. E. Humphreys, {\it Reflection groups and Coxeter
groups}, Cambridge University Press (1990).

\bibitem{hs}  T. Hausel and B. Strumfels ,  {\it Toric Hyperkahler
    Varieties},  Documenta Mathematica {\bf 7}(2002), 495-534.

\bibitem{iv} B. Iversen, {\it The Geometry of Algebraic Groups}, 
Adv. Math. {\bf 20} (1976), 57-85.





\bibitem{mer} A. S. Merkurjev, {\it Comparison between equivariant and
  ordinary $K$-theory of algebraic varieties}, Algebra i Anliz {\bf 9}
  (1997), 175-214. Translation in St. Petersburg Math. J. {\bf 9}
  (1998), 815-850.

\bibitem{st} R. Steinberg, {\it On a theorem of Pittie}, Topology
{\bf 14} (1975), 173-177.

\bibitem{str} E. Strickland, {\it Equivariant cohomology of the
wonderful group compactification},  J. Algebra, {\bf 306} (2006),
610-621.


\bibitem{su} P.Sankaran and V.Uma, {\it Cohomology of toric bundles}, {
Comment. Math. Helv.}, {\bf 78} (2003), 540-554. Errata, {\bf 79}
(2004), 840-841.

\bibitem{u} V. Uma, {\it Equivariant K-theory of compactifications of
algebraic groups}, Transformation Groups {\bf 12} (2007), 371-406.

\bibitem{u1} V.Uma, {\it Equivariant K-theory of flag varieties
    revisited and related results}, Colloq. Math. {\bf 132} (2013), 151-175.

\bibitem{vv} G.Vezzosi and A.Vistoli, {\it Higher algebraic K-theory for
actions of diagonalizable groups},  Invent. Math. {\bf 153}
(2003), 1-44.


\end{thebibliography}
\end{document}